
\input amstex
\magnification=\magstep1
\documentstyle{amsppt}
\NoRunningHeads
\vsize=8.75 true in

\parindent=.3in

\topmatter
\title
Discreteness and openness for mappings of finite distortion in the critical case $p=n-1$ 
\endtitle

\author
Enrique Villamor\\
Department of Mathematics\\
Florida International University\\
Miami, FL 33199
\endauthor

\address
Dept. of Mathematics, Florida International University, Miami, FL 
33199.
\endaddress

\email 
villamor{\@}fiu.edu
\endemail



\abstract

Let $F\in W_{loc}^{1,n}(\Omega;\Bbb R^n)$ be a mapping with non-negative Jacobian $J_F(x)=\text{det} 
DF(x)\ge 0$ a.e. in a domain $\Omega\in \Bbb R^n$. The dilatation of the mapping $F$ is defined, 
almost everywhere in $\Omega$, by the formula
$$K(x)={{|DF(x)|^n}\over {J_F(x)}}.$$
If $K(x)$ is bounded a.e., the mapping is said to be quasiregular. Quasiregular mappings are a 
generalization to higher dimensions of holomorphic mappings. The theory of higher dimensional 
quasiregular mappings began with Re\v shetnyak's theorem, stating that non constant quasiregular mappings 
are continuous, discrete and open. 

In some problems appearing in the theory of non-linear elasticity, the boundedness condition on $K(x)$ 
is too restrictive. Tipically we only know that $F$ has finite dilatation, that is, $K(x)$ is finite a.e.
 and $K(x)^p$ is integrable for some value $p$. In two dimensions, 
Iwaniec and \v Sverak [IS] have shown that $K(x)\in L^1_{loc}$ is sufficient to guarantee the conclusion of 
Re\v shetnyak's theorem.\par
For $n\ge 3$, Heinonen and Koskela [HK], showed that if the mapping is quasi-light and $K(x)\in L^p_{loc}$ for 
$p>n-1$, then the mapping $F(x)$ is continuous, discrete and open. Manfredi and Villamor [MV] proved a similar 
result without assuming that the mapping $f(x)$ was quasi-light. The result is known to be false, see [Ball], 
when $p<n-1$.\par
In this paper we attempt to improve in those results. In particular, 
we will deal with the case $p=n-1$ for $n\ge 3$, and will 
assume that our mapping $F(x)$ is quasi-light, that is, the inverse image of any point is compact in $\Omega$. 

Our approach will be different from the ones used in [MV] and [HK]. It is more geometrical in nature and uses 
the method of extremal length.

\endabstract

\endtopmatter

\document

\beginsection{$\S$1. Introduction.}

Let $\Omega\subset \Bbb R^n$, $n\ge 2$, be a domain and $F\colon \Omega\to \Bbb R^n$ be a mapping in the 
Sobolev space $W^{1,n}_{loc}(\Omega;\Bbb R^n)$ of mapping in $L^n_{loc}(\Omega;\Bbb R^n)$ whose distributional 
derivatives belong to $L^n_{loc}(\Omega;\Bbb R^n)$. The differential of $F$ at a point $x$ is denoted by 
$DF(x)$, its norm is defined by
$$\|DF(x)\|=\sup \lbrace |DF(x)\,h|\colon h\in \Bbb R^n,\,\,\,\|h\|=1\rbrace,$$
and its Jacobian determinant is defined as $J_F(x)=\text{det } DF(x)$. We assume in the rest of this paper that $F$ 
is orientation preserving, meaning that $J_F(x)\ge 0$ for a.e. $x\in\Omega$. The dilatation of $F$ at a point 
$x$ is defined by the ratio
$$K(x)={{\|DF(x)\|^n}\over {J_F(x)}}.$$
If $K(x)\in L^\infty (\Omega;\Bbb R^n)$, then $F$ is said to be a quasiregular mapping. We will say that $F$ 
is a mapping with finite dilatation if 
$$1\le K(x) <\infty$$
for almost every $x\in \Omega$ that is, except possibly for a set of measure zero in $\Omega$. We will 
follow the convention that for a.e. $x\in \Omega$ we have that if $J_F(x)=0$, then $DF(x)=0$. \par
A basic result in the theory of quasiregular mappings, states that they are discrete and open or constant. 
Vodopyanov and Goldstein [VG], proved that mappings of finite dilatation are continuous and 
have monotone components in the sende of Dirichlet. An example of J. Ball [Ball] shows that 
there are mappings satisfying that $K(x)\in L^p_{loc}$ for every $p<n-1$ that fail to be discrete.\par
A theorem of Iwaniec and \v Sverak [IS], shows that in the plane, a mapping with integrable 
dilatation, $K(x)\in L^1_{loc}$, can be expressed as the composition of an analytic mapping with an 
homeomorphism. It follows from that, that the mapping is discrete and open. The proof in [IS] is based on 
the solution of the linear two dimensional Beltrami equation and does not generalize in an obvious way to higher 
dimensions. Iwaniec and \v Sverak show that for those mapping a Stoiliv's decomposition holds. They 
conjectured in their paper that for higher dimensions, if $K\in L^{n-1}_{loc}$ then $F$ is 
discrete and open.\par
In the higher dimensional case, $n\ge 3$, Villamor and Manfredi [MV], proved that if 
$K\in L^{p}_{loc}$ for any $p>n-1$ then $F$ is discrete and open.\par
Our main goal in this paper is to prove the following result.
\proclaim{Theorem 1} Let $F\in W^{1,n}_{loc}(\Omega;\Bbb R^n)$ be a quasi-light, 
nonconstant mapping whose dilatation $K(x)$ 
is in $L^{n-1}_{loc}(\Omega)$. Then the mapping $F$ is discrete and open.
\endproclaim

The rest of the paper is constituted as follows. In section $\S$2 we will go through some rpeliminaries 
that will be needed in the rest of the paper. Section $\S$3 we will talk about some necessary conditions that 
enable us to use the change of variable formula for integrals. In Section $\S$4 we will prove the boundedness 
of the counting multiplicity function $N(F,y;\Omega)$ for the mapping $F$ and some other asymptotic 
conditions on it which will still ensure the conclusion of our main results. In Section $\S$5 we will give 
a necessary condition for the general case (i.e. not requiring that the mapping $F$ is quasi-light), which 
will allow us to prove the openess and discreteness of the mapping $F$. Section $\S$ will be devoted to show that 
the weight $w(y)=(\ln\,\ln {1\over{|y|}})^n$ is an $A_n$ weight in the sense of Muckhenhoupt. In section $\S$7 
we will prove that certain weighted variational $(n-1)$ capacity is zero for $F^{-1}(a)$, and finally in 
Section $\S$8 we will show that from that it follows that the mapping $F$ is discrete and open, which will 
end the proof of our main result.

\beginsection{ $\S$2.  Preliminaries}

{\bf 2.1 Hausd\" orff Measures.}
\vskip 0.2in
Let $s$ be a positive number and $0<\delta\le \infty$. Let $E\subset \Bbb R^n$ we define
$$\Lambda_s^\delta (E)=\inf \sum_i r_i^s,$$
where the infimum is taken over all coverings of $E$ by balls $B_i$ with radius $r_i$ not exceeding $\delta$. 
$\Lambda_s^\delta $ is an outer measurewhich in general fails to be additive on families of disjoint compact sets. 
Therefore, we define the $s$-Hausd\" orff measure of $E$ as
$$\Lambda_s(E)=\sup_{\delta>0} \Lambda_s^\delta (E)=\lim_{\delta\to 0+} \Lambda_s^\delta (E).$$
The measure $\Lambda_s$ is a Borel regular measure. That is, it is an additive measure on the Borel 
sets of $\Bbb R^n$ and for each Lebesgue measurable set $E\subset \Bbb R^n$ there is a Borel set 
$G$ such that $E\subset G$ and $\Lambda_s(G)=\Lambda_s(E)$.\par
For any set $E$, it is clear that $\Lambda_s(E)$ is a non-increasing function of $s$. Furthermore, if $s<t$, then 
$$\Lambda_s^\delta(E)\ge \delta^{s-t}\,\,\Lambda_t^\delta (E),$$
which implies that if $\Lambda_t(E)$ is positive, then $\Lambda _s(E)$ is equal to infinity. Thus, there is a 
unique value, called the Hausd\" orff dimension of the set $E$ ($\text{dim}_H (E)$), such that 
$\Lambda_s(E)=\infty$ if $0\le s<\text{dim}_H (E)$ and $\Lambda_s(E)=0$ if $\text{dim}_H (E)<s<\infty$. For 
our purposes we require however a slightly more general concept than the $s$-Hausd\" orff measure. Namely, 
let $h$ be a real valued increasing function on the interval $[0,1)$ with $\lim_{t\to 0} h(t)=0$. We define 
the $h$-Hausd\" orff measure of $E$ by 
$$\Lambda_h(E)=\sup_{\delta>0}\,\,\,inf\sum_i h(r_i),$$
where the infimum is again taken over all the coverings of the set $E$ by balls $B_i$ with radius $r_i$ not exceeding 
$\delta$. The measures $\Lambda_h$ are still Borel regular measures in $\Bbb R^n$, see [HKM] for more references. 
The choice $h(t)=t^s$ gives the $s$-Hausd\" orff measures $\Lambda _s$ defined above.
\vskip 0.2in

{\bf 2.2 Modulus.}
\vskip 0.2in
In the following, by a curve we mean a non-point locally rectifiable curve in $\Bbb R^n$. 
Let $\Gamma$ be a family of curves. We shall say that a Borel non-negative measurable function $\rho$ is 
$\Gamma$-admissible if $\int_\gamma \rho\,ds\ge 1$ for evevry $\gamma \in\Gamma$, where $s$ is the 
arc-length parameter of the curve $\gamma$. \par
Let $w$ be a non-negative measurable function in $\Bbb R^n$. We define the weighted $p$-module of $\Gamma$ by 
$$M_P^w(\Gamma)=\inf\lbrace \int_{\Bbb R^n} \rho^p\,w\,dx\colon \rho \text{ is } \Gamma-\text{ admissible}\rbrace.$$
We say that a non-negative measurable function $w$ is a weight satisfying the Muckenhoupt $A_p$ condition 
in $\Omega$ if
$$\sup_{\Bbb B\subset \Omega} \lbrace {1\over{|\Bbb B|}}\, \int _{\Bbb B} w\,dx\rbrace \,\,
\lbrace {1\over{|\Bbb B|}}\, \int _{\Bbb B} w^{{1\over{(1-p)}}}\,dx\rbrace^{p-1}<\infty,$$
where $\Bbb B$ is a ball and $|\Bbb B|$ stands for its volume. We denote it by $w\in A_p(\Omega)$. \par
Let $X$ be a bounded set in $\Bbb R^n$ containing the origin. We denote by $\Lambda({0})$ the family of rectifiable 
curves in $X$ ending at ${0}$. Then we have the following lemma.
\proclaim{Lemma 2.1} Let $w\in A_p(\Bbb R^n)$. Then the modulus $M_p^w(\Lambda ({0}))=0$ if and only if
$$\int_{|x|<1} |x|^{(1-n)p'}\,w(x)^{{1\over{(1-p)}}}\,dx =\infty,$$
where ${1\over p}+{1\over{p'}}=1.$
\endproclaim
For a proof of this Lemma, see [HH].
\vskip 0.2in
{\bf 2.2 Capacities.} A good reference for all the results in this subsection is [HKM]. We pass to define 
the weighted variational $p$-capacity of a compact set $K\subset\Omega$.
\proclaim{Definition} Let $K$ be a compact subset of $\Omega$. Define
$$W(K,\Omega)=\lbrace \mu\in C_0^\infty(\Omega)\colon \mu\ge 1 \,in\, K\rbrace.$$
We define the weighted variational $p$-capacity of $K$ as follows
$$(p,w)-\text{cap}\,(K,\Omega)= \inf_{\mu \in W(K,\Omega)}\lbrace \int_\Omega |\nabla \mu|^p\,w\, dx\rbrace.$$
\endproclaim

For any open set $U\subset \Omega$ we define
$$(p,w)-\text{cap}\,(U,\Omega)=\sup \lbrace (p,w)-\text{cap}\,(K,\Omega)\colon K\subset U \text{ compact}\rbrace$$
and finally, for an arbitrary set $E\subset\Omega$ we define
$$(p,w)-\text{cap}\,(K,\Omega)=\inf \lbrace (p,w)-\text{cap}\,(K,\Omega)\colon E\subset U\subset \Omega, \,\,U
\text{ open}\rbrace.$$
The following two properties are immediate consequences of the definition.\par
(i) If $E_1\subset E_2$, then we have that $(p,w)-\text{cap}\,(E_1,\Omega)\le (p,w)-\text{cap}\,(E_2,\Omega),$

and

(ii) If $\Omega_1\subset \Omega_2$ are open and $E\subset \Omega_1$, then we have that 
$(p,w)-\text{cap}\,(E,\Omega_2)\le (p,w)-\text{cap}\,(E,\Omega_1).$

A set $E\subset \Bbb R^n$ is said to be of variational $(p,w)$-capacity zero if 
$(p,w)-\text{cap}\,(E\cap\Omega,\Omega)=0$ for all open sets $\Omega\subset \Bbb R^n$. In this case we will write 
that $(p,w)-\text{cap}\,(E)=0$.\par
We will later need the following lemma, for reference see [HKM, p. 37].

\proclaim{Lemma 2.2} There is a positive constant $C$ independent of $p$, $w$, $x_0$ and $\tau$, such that 
$${1\over C}\,\, \Bigl( \int_{\Bbb B^n(x_0,r)} w(x)\,dx\Bigr)\,\tau^{-p}\le (p,w)-\text{cap} (\Bbb B(x_0,\tau),
\Bbb B(x_0,2\,\tau)),\tag (2.1)$$
for some positive weights, see [HKM].
\endproclaim

It is important to clarify that we could have taken the closure of the ball $\bar \Bbb B(x_0,\tau)$ instead, since 
although the $(p,w)$-$\text{cap} (G;\Omega)$ differs from $(p,w)$-$\text{cap} (\bar G;\Omega)$ when $G$ is 
an open subset of $\Omega$, when $\Bbb B$ is an open ball such that $\bar\Bbb B\subset \Omega$ we have that 
$(p,w)$-$\text{cap} (\Bbb B;\Omega)=(p,w)$-$\text{cap} (\bar\Bbb B;\Omega)$.\par
The inequality (2.1) does not hold for all weights $w$. We will show that for a particular weight $w$, 
we can follow the steps of Lemma 2.2 in [HKM] and show that for this weight the inequality in Lemma 2.2 
holds. Let us define our weight $w$ as follows,
$$w(x)= \Biggl( \ln\,\ln {1\over{|F(x)|}}\Biggr)^{n-1},$$
where we have that $|F(x)|<1$ for all $x$. 
\demo{Proof} To prove the inequality on (2.1) and according to [HKM], we 
need to show that the measure $\mu$ defined as $d\mu(x)=w(x)\,dx$, satisfies the following Poincare type 
inequality
$$\int_{\Bbb B^n(x_0,2r)} \eta (x)^{n-1}\, d\mu(x) \le C\, r^{n-1}\, 
\int_{\Bbb B^n(x_0,2r)} |\nabla\eta (x)|^{n-1}\, d\mu(x),$$
for each $\eta(x)\in C_0^\infty(\Bbb B^n (x_0,2r))$. $C$ is a constant independent of $x_0$, $r$, $\eta$ 
and $w$. In order to prove this, we will use a result proved in [MV]. Namely, 
$$V(x)=\eta(x)\,\Biggl( \ln\,\ln {1\over{|F(x)|}}\Biggr)$$
is a function in the Sobolev space $W_0^{1,n-1} (\Bbb B^n(x_0,2r))$. Thus, $V(x)$ satisfies a Poincare inequality 
with respect to the euclidean metric. Thus,
$$\int_{\Bbb B^n(x_0,2r)} \eta(x)^{n-1}\,\Bigl( \ln\,\ln {1\over{|F(x)|}}\Bigr)^{n-1} \,dx $$
$$\le C\, R^{n-1}\, \int_{\Bbb B^n(x_0,2r)} |\nabla \Biggl(\eta(x)^{n-1}\, \ln\,\ln {1\over{|F(x)|}}
\Biggr)|^{n-1} \,dx.$$
Using the product rule on the right hand side of the above inequality, we obtain that
$$\int_{\Bbb B^n(x_0,2r)} |\nabla \Biggl(\eta(x)^{n-1}\, \ln\,\ln {1\over{|F(x)|}}
\Biggr)|^{n-1}\, dx $$

$$\le 2^n\, \{ \int_{\Bbb B^n(x_0,2r)} |\nabla\eta(x)|^{n-1} 
\Biggl(\ln\,\ln {1\over{|F(x)|}}
\Biggr)^{n-1}\, dx$$

$$ + \int_{\Bbb B^n(x_0,2r)} \eta (x)^{n-1}\,|\nabla \Biggl( \ln\,\ln {1\over{|F(x)|}}
\Biggr)|^{n-1}\, dx \}.\tag (2.3)$$

By [MV], the second term on the right hand side of the above inequality is bounded by

$$\int_{\Bbb B^n(x_0,2r)} \eta (x)^{n-1}\,|\nabla \Biggl( \ln\,\ln {1\over{|F(x)|}}
\Biggr)|^{n-1}\, dx$$

$$\le C\, \Biggl(\int_{\Bbb B^n(x_0,2r)} |\nabla  \eta(x)|^n\,K(x)^{n-1}\,dx\Biggr)^{{{n-1}\over 
n}}\,\,\Biggl( \int_{\Bbb B^n(x_0,2r)} K(x)^{n-1}\,dx\Biggr)^{{1\over n}}.$$

Since by hypothesis, $K\in L^{n-1}$ and $\eta\in C_0^\infty(\Bbb B^n(x_0,2r))$, 
the right hand side of the above inequality in uniformly bounded and thus, 
among the two terms on the right hand side of (2.3), the significant one is the first term, namely with possibly 
another constant $C$, we have that 

$$\int_{\Bbb B^n(x_0,2r)} |\nabla \Biggl(\eta(x)^{n-1}\, \ln\,\ln {1\over{|F(x)|}}
\Biggr)|^{n-1}\, dx $$
$$\le C\,r^{n-1}\,  \int_{\Bbb B^n(x_0,2r)} |\nabla\eta(x)|^{n-1} 
\Biggl(\ln\,\ln {1\over{|F(x)|}}
\Biggr)^{n-1}\, dx$$

obtained the desired Poincare inequality that gives the inequality in Lemma 2.2.
\enddemo

\vskip 0.2in

To end section $\S$2, we recall that for positive weights $w$, see [HK], and for a set $E\subset \Bbb R^n$, if 
we denote by 
$$\Gamma(E)=\lbrace \gamma \text{ rectifiable curves in } \Bbb R^n \text{ ending at a point in E}
\rbrace$$ 
we have that
$$(p,w)-\text{cap} (E)=0 \text { if and only if } M_p^w(\Gamma(E))=0.$$

\beginsection{ $\S$3.  Necessary conditions for the change of variable integral formula }

We are  going to need a change of variable formula for integrals. Let $G$ be a domain in $\Bbb R^n$. A continuous 
mapping $F\colon G\to \Bbb R^n$ is said to satisfy Lusin's condition (N) if $|F(A)|=0$ whenever 
$A\subset G$ and $|A|=0$, where $|A|$ denotes the $n$-dimensional Lebesgue measure of $A$. Several necessary 
conditions have been found for a mapping $F$ to satisfy Lusin's condition (N). For our purposes we will need 
the following result on that direction which can be found in [MZ].\par
\proclaim{Theorem} Let $F$ be a mapping in the Sobolev space $W^{1,n}(\Omega;\Bbb R^n)$ which is 
continuous and such that $J_F(x)>0$ almost everywhere in $\Omega$, then $F$ satisfies condition (N) on 
$\Omega$.
\endproclaim
We will use this result in this paper, since our result is local in natureand we are dealing with mappings 
of finite dilatation which are continuous. The condition $J_F>0$ a.e. can be assumed without loss of 
generality from the sense preserving hypothesis on the mapping $F$, i.e. $J_F(x)\ge 0$ a.e. since 
otherwise $F$ will be identically constant on a set of positive measure and locally we are assuming that 
our amppings are non constant.\par
More precisely, let us suppose that $J_F\ge 0$ a.e. does not imply that $J_F>0$ a.e. Then, there exists 
an open set $U$ in $\Omega$ with $|U|>0$ such that for any $x\in U$ we have that $J_F(x)=0$. We know that 
$1\le K(x)={{|DF(x)|^n}\over {J_F(x)}}<\infty$ a.e in $\Omega$. Therefore, $|DF(x)|=0$ a.e. in $U\subset \Omega$, 
which implies that $F$ will be constant in an open subset of $U$ of positive measure.\par

The following result can be found in [MZ].
\proclaim{Theorem 2.3} Let $F\in W^{1,n}_{loc} (\Omega;\Bbb R^n)$ satisfy the Lusin's condition (N) locally 
in $\Omega$. Then we have
$$\int_{\Omega} \mu(F(x))\,J_F(x)\, dx =\int_{F(\Omega)} \mu(y)\,N(f,y;\Omega)\,dy,$$
where $N(F,y;\Omega)$ is the multiplicity function corresponding to $F$ defined as the number (possibly infinity) 
of points in $F^{-1}(y)\cap \Omega$.
\endproclaim

\beginsection{ $\S$4.  Boundedness of $N(F,y;\Omega)$}

In this section we will consider different conditions on the multiplicity function $N(F,y;\Omega)$ of the 
mapping $F$ which will lead to show that the mapping is discrete and open.\par

First, we will talk about the Brouwer degree of the mapping $F$. The Brouwer degree $d(f,\Omega,p)$ of $F$ 
with respect to $\Omega$ at a point $p\in \Bbb R^n\setminus F(\partial\Omega)$ is a well-defined 
integer depending only on the bounday values of the mapping $F$. It is well known, see [B], that if $F$ is 
smooth ($C^1$), the topological degree of a mapping $F$ $f\colon\Omega\to\Bbb R^n$ can be defined in a connected 
component $V$ of $\Bbb R^\setminus F(\partial\Omega)$ by
$$deg(F,\Omega,V)=\int_\Omega \rho(F(x))\,J_F(x)\,dx$$
where $\rho\in C_0^\infty(\Bbb R^n)$, such that its support is contained in $V$ and satisfying the condition 
$\int_{\Bbb R^n} \rho(y)\,dy=1$. Our next goal is to show that this is still true for mappings $F$ continuous in 
$\bar\Omega$ and belonging to the Sobolev space $W_{loc}^{1,n}(\Omega)$.

\demo{Proof} Let $V$ be a connected component of $\Bbb R^n\setminus F(\partial\Omega)$.  It is clear that
$d(F,\Omega,p)$ is independent of $p\in V$.  We denote this common values for all $p\in V$ by $d(F,\Omega,V)$.  Without
loss of generality we can assume that $p=(0,0,\ldots,0)$.  Let $A$ be a subset of $\Omega$ containing 0 with
$A\subset\bar A\subset \Omega$, such that the support of $rho$ lies in a connected component of $\Bbb R^n\setminus
F(\bar\Omega\setminus A)$.  Now we smooth tha mapping $F$ on the set $A$.  Let $\epsilon
<\text{dist}(A,\partial\Omega)$ and define the mollifiers, $g_\epsilon(x)= {1\over{\epsilon^n}}\,g({x\over\epsilon})$
for any $x\in \Bbb R^n$ and $g\in C_0^\infty(\Bbb R^n)$ such that the suport of $g$ is contained in the unit ball $\Bbb
B^n(0,1)$ of $\Bbb R^n$ and $\int{\Bbb R^n} g(x)\, dx=1$.  let us define $F_\epsilon (x)$ as the convolution of the
function $g_\epsilon$ with the mapping $F$ as follows $$F_\epsilon(x)= \int_{\Bbb R^n} g_\epsilon (x-y)\,F(y)\,dy,$$
for any $x\in A$.  Since the convolution of two functions satisfy the commutative property, we obtain that
$$F_\epsilon(x)= \int_{\Bbb R^n} F(x-y)\,{1\over{\epsilon^n}}\, g({x\over\epsilon})\,dy.$$ By performing the change of
variable $\tilde y={y\over{\epsilon}}$ we obtain the following formula for $F_\epsilon$, $$f_\epsilon(x)= \int_{\Bbb
R^n} F(x-\epsilon \tilde y)\,g(\tilde y)\, d\tilde y.$$ Thus, $F_\epsilon$, the convolution of the mapping $F$ with the
function $g_\epsilon$, belongs to $C(\bar A)\cap W^{1,n}(A)$.  We will show that $F_\epsilon$ converges uniformly on
compact subsets of $A$ to the mapping $F$.  For this, observe that $$|F_\epsilon(x)-F(x)|=|\int_{\Bbb R^n} F(x+\epsilon
y)\,g(y)\,dy -\int_{\Bbb R^n} F(x)\,g(y)\,dy|,$$ where to simplify the notation, we have replaced $-\epsilon$ by
$\epsilon$ and $\tilde y$ by $y$.  Hence, we have that $$|F_\epsilon(x)-F(x)|=|\int_{\Bbb R^n} (F(x+\epsilon
y)-F(x))\,g(y)\,dy|,$$ obtaining that $$|F_\epsilon(x)-F(x)|\le \int_{\Bbb R^n} |(F(x+\epsilon y)-F(x))|\,|g(y)|\,dy.$$
Let us denote by $G_\epsilon(x)=F(x+\epsilon\, y)$ for $x\in \Omega$.  Then, it is trivial to see that the family of
mappings $\lbrace G_\epsilon\rbrace_\epsilon$ converges uniformly to the mapping $F$ in $\Omega$ as $\epsilon$
approaches 0.  Thus, letting $\epsilon\to 0$, $|F_\epsilon(x)-F(x)|\to 0$ uniformly, which implies that $F_\epsilon$
converges to $F$ uniformly on compact subsets of $A$ as $\epsilon \to 0$.  It is also clear that for $\epsilon$ small
enough, the support of $\rho$ lies in a connected component of $\Bbb R^n \setminus F_\epsilon(\partial A)$.  Thus,
since $F_\epsilon$ is a smooth mapping, we have the follwoing formula to compute its Browwer degree,
$$\text{deg}(F_\epsilon,\Omega,0)=\int_{\Bbb R^n} \rho(F_\epsilon(x))\,J_{F_\epsilon}(x)\, dx,$$ letting $\epsilon\to
0$ we obtain the desired formula, $$\text{deg}(F,\Omega,0)=\int_{\Bbb R^n} \rho(F(x))\,J_F(x)\, dx,\tag (2.4)$$ making
use of the fact that $\rho(F_\epsilon(x))\,J_{F_\epsilon}(x)\to \rho(F(x))\,J_F(x)$ uniformly as $\epsilon \to 0$.\par
We have already established the validity in our case, of the change of variable formula, $$\int_\Omega
\mu(F(x))\,J_F(x)\,dx=\int_{F(\Omega)} \mu(y)\, N(F,\Omega,y)\,dy.$$ If we let in the above formula, $\mu=1$, we have
that $$\int_\Omega J_F(x)\,dx=\int_{F(\Omega)} \, N(F,\Omega,y)\,dy.\tag (2.5)$$ Let us construct now a family of
functions $\rho_r$ as follows.  Let $V$ be a connected component of $\Bbb R^n \setminus F(\partial \Omega)$ in
$F(\Omega)$.  Let $Theta_r$ be a continuous function in $\Bbb R^n$ satisfying that $\Theta_r(p)=1$ for all $p\in V$
such that $\text{dist}(p,\partial V) \ge {1\over r}$, with its support in $V$ and $0\le \Theta_r(p)\le 1$ for any $p\in
\Bbb R^n$.\par Let us define the functions $\rho_r$ as follows, $$\rho_r(p)={{\Theta_r(p)}\over{\int_{\Bbb R^n}
\Theta_r(p)\,dp}}.  $$ Then we have that, $$\text{deg}(F,\Omega,0)=\int_{\Bbb R^n} \rho_r(F(x))\,J_F(x)\,dx.$$ By the
definition of $\rho_r$ we have that $$\text{deg}(F,\Omega,0)={1\over{ \int_{\Bbb R^n} \Theta_r(p)\,dp}}\, \int_{\Bbb
R^n} \Theta_r(F(x))\,J_F(x)\,dx,$$ letting $r\to\infty$ in the above equality, we obtain that
$$\text{deg}(F,\Omega,0)={1\over{ |V|}}\, \int_{F^{-1}(V)} J_F(x)\,dx.$$ Thus, we have that
$$|V|\,\text{deg}(F,\Omega,0)=\int_{F^{-1}(V)} J_F(x)\,dx.\tag (2.6)$$

Combining (2.5) amd (2.6) we obtain that

$$|V|\,\text{deg}(F,\Omega,0)=\int_{V} \, N(F,\Omega,y)\,dy,$$
which proves that a.e. in $y\in V$ we have that
$$N(F,\Omega,y)=\text{deg}(F,\Omega,0).$$

\vskip 0.2in

Let us consider an special case. Let us suppose that for our mapping $F\colon\Omega\to \Bbb R^n$ there 
exists another mapping $F_0\colon \Omega\to \Bbb R^n$ which is continuous in $\bar \Omega$ and one to one in $\Omega$. 
Let us further assume that $F|_{\partial\Omega}=F_)|_{\partial\Omega}$. Then, it is not difficult to show 
that $N(F,\Omega,y)=1$ a.e. in $F(\Omega)$, see [Ball].\par
 It is important to either show that under our hypotheis of the mapping $F$ being quasi-light, then the 
 multiplicity function $N(F,\Omega,y)$ is bounded or to impose restrictions on this multiplicity 
 function that assure us that the openness and discreteness conclusions of our results still 
 hold. Let us consider two different scenarios that allow us to conclude that our multiplicity function is 
 bounded.
 \vskip 0.2in
 {\bf 1. If the mapping $F$ under consideration is quasi-light, then $N(F,\Omega,y)$ is bounded (or 
 locally bounded in $F(\Omega)$, which for our purposes is enough, since our result is local in nature).}\par
 
 {\bf 2. By [FG], we have that if the mapping $F\in W^{1,n}(\Bbb B^n(x,R))$ and $J_F(x)>0$ a.e. then $F$ 
 has a differential almost everywhere (in the classical sense) and there exists a $R_0>0$ such that for every 
 $0<r<R_0$ we have that
 $$\text{deg} (F,\Bbb B^n(x_0,r),y)=1$$
 for every $y\in C_r$ where $C_r$ is the connected component of $\Bbb R^n\setminus F(\partial \Bbb B^n(x_0,r))$ 
 containing $y_0=F(x_0)$, see Lemma 6.5 in [FG].}
 \enddemo

\beginsection{ $\S$5.  A related result for the general case}

Since our result is local, it will be enough to concentrate in a neighborhood of the origin 
and then by means of a linear transformation the result will remain valid for any $b\in \Bbb R^n$. 
Let $\Omega$ be an open ball such that $F^{-1}(0)=E\subset \Omega$. Let $\Delta$ be the family of 
rectifiable curves in $F(\Omega)$ ending at the origin. Similarly, let $\Delta_1$ be the family of rectifiable curves 
ending at a point in $E$ and completely contained in $\Omega$. Let $\Gamma(\Delta)$ be the family of admissible 
metrics for $\Delta$, and similarly, let $\Gamma(\Delta_1)$ be the family of admissible metrics for $\Delta_1$. 

Let $\gamma\in \Delta_1$ be such that its components $\gamma(t)=\lbrace (x_1(t),x_2(t),\ldots,x_n(t))\colon 
t\in [a,b]\rbrace$, then it follows that $F(\gamma)\in \Delta_1$ and 
$$F(\gamma(t))
=\lbrace (y_1(t),y_2(t),\ldots,y_n(t))\colon y_i(t)=F_i(x_1(t),x_2(t),\ldots,x_n(t))$$
$$, i=1,2, \ldots,n 
\text{ and }t\in [a,b]\rbrace.$$ 
Let 
$$d\bar s(y)=\sqrt{\Bigl( {{\partial y_1}\over{\partial t}}\Bigr)^2+
\ldots+\Bigl( {{\partial y_n}\over{\partial t}}\Bigr)^2}.$$ 
By the definition of $\gamma_i,\,i=1,2,\ldots,n$ and using 
the chain rule for differentiation we have that
$${{\partial y_i}\over{\partial t}}=\sum_{j=1}^n {{\partial F_i}\over{\partial x_j}}
\,{{\partial x_j}\over{\partial t}},$$
for any $i=1,2,\ldots,n$. Thus, we obtain that
$$d\bar s(y)=\sqrt{\Biggl[\Bigl( \sum_{i=1}^n{{\partial F_1}\over{\partial x_i}}\,\dot x_i\Bigr)^2+
\ldots+\Bigl(\sum_{i=1}^n {{\partial F_n}\over{\partial x_i}}\,\dot x_i\Bigr)^2\Biggr]},$$
where $\dot x_i$ stands for ${{d x_i}\over {dt}}$.Applying Schwartz's inequality we obtain that
$$d\bar s(y)\le\sqrt{\Biggl[ \sum_{i=1}^n\Bigl({{\partial F_1}\over{\partial x_i}}\Bigr)^2\,\sum_{i=1}^n 
(\dot x_i)^2+
\ldots+\sum_{i=1}^n\Bigl({{\partial F_n}\over{\partial x_i}}\Bigr)^2\,\sum_{i=1}^n 
(\dot x_i)^2\Biggr]}.$$
Using the fact that $ds(x)=\sqrt{\sum_{i=1}^n (\dot x_i)^2}$ we obtain that
$$d\bar s(y)\le\sqrt{\Biggl[ \sum_{i=1}^n\Bigl({{\partial F_1}\over{\partial x_i}}\Bigr)^2\,(ds(x))^2+
\ldots+\sum_{i=1}^n\Bigl({{\partial F_n}\over{\partial x_i}}\Bigr)^2\,(ds(x))^2\Biggr]}$$

$$=ds(x)\,\sqrt{\Biggl[ \sum_{i=1}^n\Bigl({{\partial F_1}\over{\partial x_i}}\Bigr)^2+
\ldots+\sum_{i=1}^n\Bigl({{\partial F_n}\over{\partial x_i}}\Bigr)^2\Biggr]}.$$
We also have that,
$$\sum_{i=1}^n\Bigl({{\partial F_j}\over{\partial x_i}}\Bigr)^2\le |DF(x)|^2$$
for every $j=1,2,\ldots,n$, which implies that
$$d\bar s(y)\le \, \sqrt {n\, |DF(x)|^2}\,ds(x)= \sqrt n\, |DF(x)|\,ds(x).$$
We will use this inequality to show that if $\rho$ is an admissible metric for the family of curves $\Delta$ 
then the metric defined ny $\sqrt n\, (\rho\circ F)(x)\,|DF(x)|$ is an admissible metric for the family of curves 
$\Delta_1$. For this, let $\gamma\in \Delta_1$ then we have that $\hat\gamma=F(\gamma)\in \Delta$. Since $\rho$ 
is an admissible metric for $\delta$, we have that $1\le \int_{\hat\gamma} \rho(y)\,d\bar s(y)$ which implies that 
$$1\le \sqrt n \, \int_\gamma \rho(F(x))\, |DF(x)|\,ds(x).$$
Since $\gamma$ was any arbitrary curve in $\Delta_1$, we have that the metric $\sqrt n\, \rho(F(x))\,|DF(x)|$ 
is admissible for the family $\Delta_1$, that is, we have concluded that 
$$\sqrt n\, \rho(F(x))\,|DF(x)|\in 
\Gamma(\Delta_1)$$ 
whenever $\rho\in \Gamma(\Delta)$.\par
Let us recall the definition of weighted module of order $p$ when the weight function is identically equal to 1. 
We will denote by $M_p$. It follows from the definition that
$$M_p(\Delta_1)\le \sqrt{n^p}\,\int_\Omega (\rho(F(x)))^p\,|DF(x)|^p\,dx.$$
Multiplying and dividing the integrand on the right hand side of the above inequality by $K(x)^{{p\over n}}$ 
we obtain that
$$M_p(\Delta_1)\le \sqrt{n^p}\,\int_\Omega (\rho(F(x)))^p\,|DF(x)|^p\,K(x)^{{p\over n}}\, 
{1\over{K(x)^{{p\over n}}}}\,dx,$$
applying H\"older's inequality we have that
$$M_p(\Delta_1)\le C\,\Biggl(\int_\Omega (\rho(F(x)))^n\,{{|DF(x)|^n}\over{K(x)}}\,dx\Biggr)^{{p\over n}}\,\,
\Biggl( \int_\Omega\Bigl(K(x)^{{p\over n}}\Bigr)^{{n\over{n-p}}}\,dx\Biggr)^{{{n-p}\over n}},$$
where $C$ is a constant depending only on $n$ and $p$. Now, since $J_F(x)={{|DF(x)|^n}\over {K(x)}}$ we have that
$$M_p(\Delta_1)\le C\,\Biggl(\int_\Omega (\rho(F(x)))^n\,J_F\,dx\Biggr)^{{p\over n}}\,\,
\Biggl( \int_\Omega K(x)^{{p\over {n-p}}}\,dx\Biggr)^{{{n-p}\over n}}.\tag (5.2)$$
We are going to use the change of variable formula in the first integral on the right hand side of the 
above inequality. Thus, we have
$$\int_\Omega (\rho(F(x)))^n\,J_F\,dx= \int_{F(\Omega)} (\rho(y))^n\,N(F,\Omega,y)\,dy.$$
Now we will study different conditions on the growth of the multiplicity function $N(F,\Omega,y)$ which 
will still guarantee that by choosing a convenient sequence of admissible metrics $\rho_\eta$ in the above 
equality and then taking the limit, that right hand side goes to zero. From now on, $C$ will denote possibly 
different constants independent of $F$, $\Omega$ and $\rho$. Summarizing, we have that 
$$M_p(\Delta_1)\le C\,\Biggl(\int_{F(\Omega)} (\rho(y))^n\,N(F,\Omega,y)\,dy\Biggr)^{{p\over n}}\,\,
\Biggl( \int_\Omega K(x)^{{p\over {n-p}}}\,dx\Biggr)^{{{n-p}\over n}}.\tag (5.3)$$
AT this point, let us assume that $N(F,\Omega,y)$ is bounded. It is also clear that we can replace in (5.3) 
$N(F,\Omega,y)$ by $N(F,\Omega,r)=\sup_{y\in \partial \Bbb B^n(0,r)} N(F,\Omega,y)$, which is now a radial function 
of $r$. That is, we have the alternating formula to (5.3),
$$M_p(\Delta_1)\le C\,\Biggl(\int_0^R\,\int_{\partial \Bbb B^n(0,r)} (\rho(y))^n
\,N(F,\Omega,r)\,dS(r)\,dr\Biggr)^{{p\over n}}\,\,
\Biggl( \int_\Omega K(x)^{{p\over {n-p}}}\,dx\Biggr)^{{{n-p}\over n}}.\tag (5.4)$$

\vskip 0.2in
{\bf Case 1. $N(F,\Omega,y)$ is bounded.}
\vskip 0.2in
Then by 5.3 we have that,
$$M_p(\Delta_1)\le C\,\Biggl(\int_{F(\Omega)} (\rho(y))^n\,dy\Biggr)^{{p\over n}}\,\,
\Biggl( \int_\Omega K(x)^{{p\over {n-p}}}\,dx\Biggr)^{{{n-p}\over n}}.$$
Let $p=n-1+\epsilon>n-1$ and observe that
$${p\over{n-p}}={{n-1+\epsilon}\over{1-\epsilon}}=n-1+{{n\epsilon }\over{1-\epsilon}}=n-1+\delta,$$
where $\delta={{n\epsilon}\over{1-\epsilon}}$. hence we have that
$$M_{n-1+\epsilon}(\Delta_1)\le C\,\Biggl(\int_{F(\Omega)} (\rho(y))^n\,dy\Biggr)^{{{n-1+\epsilon}\over n}}\,\,
\Biggl( \int_\Omega K(x)^{n-1+\delta}\,dx\Biggr)^{{{1-\epsilon}\over n}}.\tag (2.6)$$
Taking infimums in the above inequality we obtain that
$$M_{n-1+\epsilon}(\Delta_1)\le C\,(M_n(\Delta))^{{{n-1+\epsilon}\over{n}}}\,
\Biggl( \int_\Omega K(x)^{n-1+\delta}\,dx\Biggr)^{{{1-\epsilon}\over n}}.\tag (2.7)$$
Next, we are going to show that $M_n(\Delta)$ is equal to zero. For this, we 
will use Lemma II.2.1  with the weight $w$ identically equal to one, whic obviously belongs to the 
Muckenphout class $A_n$. The Lemma states that $M_n(\Delta)=0$ if and only if 
$$\int_{|x|<1} |x|^{(1-n){n\over{n-1}}}\,dx=\infty.$$
Using spherical coordinates in $\Bbb R^n$ it is trivial to show that the above integral becomes a 
divergent improper integral and we are done.\par
We would like to remind here, see [MV], that the reason to have an $\epsilon$ in our 
argument is because a classical result states that if $M_p(\delta_1)=0$ then the Hausdorff diemnsion of $E$ 
is less than or equal to $n-p$. In our situation, we want to ensure that $M_{n-1+\epsilon} (\Delta_1)=0$, thus 
according to this classical result, $\text{Hausdorff dimension } (E)\le 1-\epsilon<1$ which 
implies the discreteness of the set $E$ and now by Titus and Young [TY], the openness of $F$. \par
AT this point, we can not conclude yet that $M_{n-1+\epsilon} (\Delta_1)=0$ since we only know that 
$K(x)\in L^{n-1}_{loc} (\Omega)$ and the integral that appears on the right hand side of (2.7) is 
$\int_\Omega K(x)^{n-1+\delta}\,dx$. Since we can assume without loss of generality since our results 
are local, that $L^{n-1+\delta}_{loc} (\Omega)\subset L^{n-1}_{loc} (\Omega)$. At this point, and since the case 
$K(x)\in L^{n-1+\delta}_{loc} (\Omega)$ was completely settled by [MV], we can assume that 
$K(x)\in L^{n-1}_{loc} (\Omega)\setminus L^{n-1+\delta}_{loc} (\Omega)$. Thus, we are facing on the right 
hand side of (2.7) a product of the form zero times infinity. Our goal now will be to find which condition 
is needed for that indeterminate product to be zero. For this, let us define the following metric, 
$\rho=|\nabla \Bigl(\ln {1\over{|y|}}\Bigr)^\delta|$ with $\delta$ positive and strictly less 
than $1-{1\over n}$. 

We want to show that, somehow we can 
use these metrics to obtain an estimate of the modulus $M_n(\Delta)$. For this, we observe that $\Delta$ is 
the union of $\bigcup_{\eta>0}\Delta_\eta$, where $\Delta_\eta $ denotes the family of rectifiable curves
in $\Bbb B^n(0,\eta)$ joining the origin with a point on the boundary of $\Bbb B^n(0,\eta)$. It is immediate 
to show that the modulus of this family of curves is bigger than the modulus of the family of curves in $\Bbb R^n$ 
with one end point at the origin and another at a boundary point of $\Bbb B^n(0,\eta)$. And this as we will 
see later in our argument will be enough for our arguments, thus we can restrict ourselves to estimate 
the modulus of the former family.\par
It is also clear that for each $\gamma\in \Gamma(\Delta)$ we have that $\int_\gamma \rho\,ds=\infty \ge 1$. Thus, 
it is admissible for $M_n(\Delta)$ and thus so are $\rho_\epsilon=\epsilon\, 
|\nabla \Bigl(\ln {1\over{|y|}}\Bigr)^\delta|$. Let us compute now $\int_{\Bbb B^n(0,1/2)} \rho_\epsilon^n(x)\,dx$. 
A straightforward computation gives us that 
$$|\nabla \Bigl(\ln {1\over{|y|}}\Bigr)^\delta|=\delta\, \Bigl(\ln {1\over{|y|}}\Bigr)^{\delta-1}\,
{1\over {|y|}}.$$
Let us denote $|y|=r$, thus we have after passing to spherical coordinates in $\Bbb R^n$
$$\int_{\Bbb B^n(0,1/2)} \rho_\epsilon^n(x)\,dx= C\,\epsilon^n\,\delta^n\,\int_0^1/2
\Bigl(\ln {1\over{r}}\Bigr)^{n(\delta-1)}\,
{1\over {r}}^n\,r^{n-1}\,dr,$$
where $C$ is a constant independent of $\delta$, $r$, and $\epsilon$. Using the change of variable 
$u=\ln {1\over{r}}$ the above integral is transformed to
$$\int_{\Bbb B^n(0,1/2)} \rho_\epsilon^n(x)\,dx= C\,\epsilon^n\,\delta^n\,\int_{\ln 2}^\infty
\Bigl(u\Bigr)^{n(\delta-1)}\,du,$$
the improper integral above converges by our choice of $\delta$, and thus we have that
$$\int_{\Bbb B^n(0,1/2)} \rho_\epsilon^n(x)\,dx= C\,\epsilon^n,$$
letting $\epsilon\to 0$ we have shown that $M_n(\Delta)=0$.\par
This provides a proof of Heinonen and Koskela's result, see [HK] in the Archive for Rational Mechanics, that 
a quasi-light mapping $F\in W^{1,n}_{loc}(\Omega;\Bbb R^n)$ with dilatation $K(x)\in L^p_{loc}(\Omega)$ 
for some $p>n-1$ is discrete and open.

\vskip 0.6in

{\bf $\S$6. The weight $w(y)=(\ln\,\ln {1\over{|y|}})^n$ is an $A_n$ Muckhenhoupt weight}

In order to prove our main result, which improves our previous ones, we need to introduce weighted 
modulus and their corresponding weighted variational capacities, and show that the weight 
$$w(y)=(\ln\,\ln {1\over{|y|}})^n$$
defined in $0<|y|<1$, satisfies the Muckenhoupt $A_n$ condition in $\Omega$.\par
We first observe that without loss of generality we can assume that $F(\Omega)\subset \Bbb B^n(0,1)$, with 
$0\in F(\Omega)$, so that $w$ is defined in $F(\Omega)$. It is also enough on the weight condition, to 
take the supremum over balls centered at 0 and show that
$$\sup_{0<r_0<1} \Biggl[{1\over {r_0^n}}\, \int_0^{r_0} \Bigl(\ln\,\ln {1\over r}\Bigr)^n\,r^{n-1}\,dr\Biggr]\,
 \Biggl[{1\over {r_0^n}}\, \int_0^{r_0} \Bigl(\ln\,\ln {1\over r}\Bigr)^{{n\over {1-n}}}\,r^{n-1}\,dr\Biggr]^{n-1},$$
 is finite. In order to show that, let us start finding a bound for the second factor on the right hand side 
 of the above expression. For simplicity, we will denote that right hand side by $A(r_0)$. Since $r\le r_0$, 
 we obtain the following inequality
 $$\ln\,\ln {1\over {r_0}}\le \ln\,\ln {1\over {r}},$$
 exponentiating, using the negative exponent ${n\over{1-n}}$, we obtain that
 $$\Biggl(\ln\,\ln {1\over {r_0}}\Biggr)^{{n\over{1-n}}}\le \Biggl(\ln\,\ln {1\over {r}}\Biggr)^{
 {n\over{1-n}}}.$$
 Thus, we have that 
 $$A(r_0)=\Biggl[{1\over {r_0^n}}\, \int_0^{r_0} \Bigl(\ln\,\ln {1\over r}\Bigr)^{{n\over {1-n}}}
 \,r^{n-1}\,dr\Biggr]^{n-1}$$
 $$\le \Biggl[{1\over {r_0^n}}\,
  \int_0^{r_0} \Bigl(\ln\,\ln {1\over {r_0}}\Bigr)^{{n\over {1-n}}}\,r_0^{n-1}\,dr\Biggr]^{n-1},$$
  since $\int_0^{r_0} dr=r_0$ we obtain that,
 $$A(r_0)\le  \Biggl[{1\over{r_0^n}}\,
  \Bigl(\ln\,\ln {1\over {r_0}}\Bigr)^{{n\over {1-n}}}\,r_0^{n}\,dr\Biggr]^{n-1}=
  \ln {1\over {r_0}}\Bigr)^{-n}.$$
We shall now compute the integral on the first factor on the right hand side of the $A_n$ condition. Namely,
$$\int_0^{r_0} \Bigl(\ln\,\ln {1\over r}\Bigr)^n\,r^{n-1}\,dr.$$
For this, we will use the following substitution, $u={1\over r}$, then we have that $r={1\over u}$ and 
$dr=-{1\over{u^2}}\,du$. Substituting this in the above integral we obtain that
$$\int_{{1\over {r_0}}}^{\infty} \Bigl(\ln\,\ln u\Bigr)^n\,{1\over {u^{n+1}}}\,du.$$
We will perform another substitution, $T=\ln u$, then $u=e^T$ and $du=e^T\,dT$, and hence the above integral becomes 
$$\int_{\ln{1\over {r_0}}}^{\infty} \Bigl(\ln T\Bigr)^n\,{1\over {e^{T(n+1)}}}\,e^T\,dT=
\int_{\ln{1\over {r_0}}}^{\infty} \Bigl(\ln T\Bigr)^n\,{1\over {e^{Tn}}}\,dT.$$
We will compute the last improper integral using the method of integration by parts. Let $u=\Bigl(\ln T\Bigr)^n$ 
and $dv=e^{-nT}\,dT$. Then, $du=n\,\Bigl(\ln T\Bigr)^{n-1}\,{1\over T}\, dT$ and $v=-{{e^{-nT}}\over n}$. Thus, 
integrating by parts, we obtain that the above integral is equal to
$$\int_{\ln{1\over {r_0}}}^{\infty} \Bigl(\ln T\Bigr)^n\,{1\over {e^{Tn}}}\,dT=
\Biggl[-\Bigl(\ln T\Bigr)^n\, {{e^{-nT}}\over n}\Biggr]^\infty_{\ln{1\over {r_0}}}$$
$$+\int_{\ln{1\over {r_0}}}^{\infty} 
-{e^{-nT}}\,\Bigl(\ln T\Bigr)^{n-1}\,{1\over T}\, dT.$$
We observe now that the second term on the right hand side of the above equality is majorized by the 
first term , so we can disregard that term. thus, we obtain that
$$\int_{\ln{1\over {r_0}}}^{\infty} \Bigl(\ln T\Bigr)^n\,{1\over {e^{Tn}}}\,dT \approx 
{{r_0^n}\over n}\, \Biggl(\ln\,\ln {1\over {r_0}}\Biggr)^n.$$
Putting both estimates together in the $A_n$ condition we obtain that
$$\sup_{0<r_0<1} \Biggl[{1\over {r_0^n}}\, \int_0^{r_0} \Bigl(\ln\,\ln {1\over r}\Bigr)^n\,r^{n-1}\,dr\Biggr]\,
 \Biggl[{1\over {r_0^n}}\, \int_0^{r_0} \Bigl(\ln\,\ln {1\over r}\Bigr)^{{n\over {1-n}}}\,r^{n-1}\,dr\Biggr]^{n-1}$$
 $$\le \sup \Biggl[ {1\over {r_0^n}}\,{{r_0^n}\over n}\,\Biggl(\ln\,\ln {1\over {r_0}}\Biggr)^n\,\,
 \Biggl(\ln\,\ln {1\over {r_0}}\Biggr)^{-n}\Biggr]={1\over n}<\infty,$$
 and this concludes our proof that our weight $w(y)=(\ln\,\ln {1\over{|y|}})^n$ belongs to the Muckenhoupt 
 class $A_n(\Omega)$.

\beginsection{ $\S$7.  Weighted $(n-1)$ variational capacities}

In this section we will move from unweighted modulus and variational capacities to the weighted versions 
of them. This will allow us to improve on our results on the openness and discreteness of mappings 
$F$ with integrable dilatation. We will also use $h$-Hausdorff-measures in our arguments in this section and 
the next.\par
We shall start by recalling the fact that if $\rho$ is an admissible metric for the family of curves $\Delta$, 
then $\sqrt n\, \rho(F(x))\, |DF(x)|$ is an admissible metric for the family of curves $\Delta_1$ as 
it was proved in section $\S$5. Let $w_1(x)$ be a positive weight defined by
$$w_1(x)=\Biggl(\ln\,\ln {1\over{|F(x|}}\Biggr)^{n-1}$$
in $\Omega$. Since our result is local we can assume without loss of generality that $|F(x)|<1$ in $\Omega$. Then 
by the definition of weighted modulus we have that
$$M_p^{w_1}(\Delta_1)\le \int_\Omega \sqrt n^p\,(\rho(F(x)))^p\,\Biggl(\ln\,\ln {1\over{|F(x|}}\Biggr)^{n-1}\,
|DF(x)|^p\,dx,$$
multiplying and dividing the integrand on the right hand side of the above inequality by 
$K(x)^{{p\over n}}$, we have that
$$M_p^{w_1}(\Delta_1)\le \int_\Omega \sqrt n^p\,(\rho(F(x)))^p\,\Biggl(\ln\,\ln {1\over{|F(x|}}\Biggr)^{n-1}\,
|DF(x)|^p\,K(x)^{{p\over n}}\, {1\over{K(x)^{{p\over n}}}}\,dx.$$
Applying H\"older's inequality, we obtain that
$$M_p^{w_1}(\Delta_1)\le \sqrt n^p\,\Biggl[\int_\Omega (\rho(F(x)))^n\,
\Biggl(\ln\,\ln {1\over{|F(x|}}\Biggr)^{{{(n-1) n}\over p}}\,
{{|DF(x)|^n}\over {K(x)}}\,dx\Biggr]^{{p\over n}}$$
$$\Biggl[ \int_\Omega K(x)^{{p\over {n-p}}}\,dx\Biggr]^{{{n-p}\over n}}.$$
Since ${{|DF(x)|^n}\over {K(x)}}=J_F(x)>0$ a.e. and using the formula for the change of variables on the first 
integral on the right hand side of the above inequality, we obtain that
$$M_p^{w_1}(\Delta_1)\le \sqrt n^p\,\Biggl[\int_{F(\Omega)} (\rho(y))^n\,
\Biggl(\ln\,\ln {1\over{|y|}}\Biggr)^{{{(n-1) n}\over p}}\,
N(F,\Omega,y)\,dy\Biggr]^{{p\over n}}$$
$$\Biggl[ \int_\Omega K(x)^{{p\over {n-p}}}\,dx\Biggr]^{{{n-p}\over n}},$$
let $p=n-1$ in the above inequality to obtain
$$M_{n-1}^{w_1}(\Delta_1)\le \sqrt n^{n-1}\,\Biggl[\int_{F(\Omega)} (\rho(y))^n\,
\Biggl(\ln\,\ln {1\over{|y|}}\Biggr)^{n}\,
N(F,\Omega,y)\,dy\Biggr]^{{{n-1}\over n}}$$
$$\Biggl[ \int_\Omega K(x)^{n-1}\,dx\Biggr]^{{1\over n}}.$$

Next, we are going to examine the integral 
$$\int_{F(\Omega)} (\rho(y))^n\,
\Biggl(\ln\,\ln {1\over{|y|}}\Biggr)^{n}\,
N(F,\Omega,y)\,dy.$$
 First, we are going to assume that our metrics $\rho$ are radial, that is, $\rho(y)=\rho (r)$ where $|y|=r$, 
 and thus by taking $N(F,\Omega,r) =\sup_{y\in \partial \Bbb B^n(0,r)} N(F,\Omega,y)$, the above integral 
 is less than or equal to
 $$C\,\int_0^{1\over 2} (\rho(r))^n\,
\Biggl(\ln\,\ln {1\over{|y|}}\Biggr)^{n}\,
N(F,\Omega,r)\,dr.$$

\vskip 0.2in
We will consider first the case $K(x)\in L^{n-1}_{loc}(\Omega)$ and $F$ quasi-light, that is the 
multiplicity function $N(F,\Omega,y)$ is essentially bounded in $\Omega$. \par
In [Hencl and Maly], it is shown that if the mapping $F\in W^{1,p}_{loc} (\Omega;\Bbb R^n)$, with $p>n-1$, be a 
continuous mapping with finite distortion which satisfies that $J_F\in L^1_{loc}(\Omega)$ and the equality 
$$\text{div}((\Psi\circ F)\,\text{ adj}\,DF)=((\text{div }\Psi)\circ F)\,J_F$$
holds in the sense of distributions in $\Omega$ for each $C^1$-vector field $\Psi$ on $\Bbb R^n$. Then for any 
$\Omega '$ relatively compact subset of $\Omega$ we have that
$$N(F,\Omega ',y)=\text{deg}(F,\Omega ', y)$$
for a.e. $y\in \Bbb R^n\setminus F(\partial\Omega ')$.\par
Observe that our mapping $F$ satisfies these hypothesis. Then taking into consideration our 
Remark 2 at the end of section $\S$4, we have that for $F$ quasi-light, choosing a point $x_0\in \Omega$ such 
that $F(x_0)=0$ we can find an open set $\Omega '$ compactly contained in $\Omega$ including the connected 
component of $F^{-1}(0)$ containing $x_0$, such that $0\notin F(\partial \Omega ')$. Let $\rho>0$ such that
$$\bar \Bbb B^n((0,\rho)\cap F(\partial \Omega ')=\emptyset$$
and define $\Omega ''$ as the connected component of $\Omega '\cap F(\Bbb B^n((0,\rho))$ which contains $x_0$. 
Then
$$N(F,\Omega '',y)=\text{deg}(F,\Omega '', y)=\text{deg}(F,\Omega '', 0)$$
if $y\in \Bbb B^n((0,\rho)$ and $N(F,\Omega '',y)=\text{deg}(F,\Omega '', y)=0$ if $y\notin \Bbb B^n((0,\rho)$. 
This shows that $N(F,\Omega '',\,\,\,)$ is essentially bounded in a neighborhood of $0$.

This implies that we have
$$M_{n-1}^{w_1}(\Delta_1)\le C\,\Biggl[\int_{F(\Omega)} (\rho(y))^n\,
\Biggl(\ln\,\ln {1\over{|y|}}\Biggr)^{n}\,dy\Biggr]^{{{n-1}\over n}}$$
$$\Biggl[ \int_\Omega K(x)^{n-1}\,dx\Biggr]^{{1\over n}}.$$
Taking infimums over all the admissible metrics $\rho$ for the family of curves $\Delta$ we have that
$$M_{n-1}^{w_1}(\Delta_1)\le \sqrt n^{n-1}\, (M_n^w(\Delta))^{{{n-1}\over n}}\,\,
\Biggl[ \int_\Omega K(x)^{n-1}\,dx\Biggr]^{{1\over n}},$$
where the positive weight $w$ is defined by 
$$w(y)=\Biggl(\ln\,\ln {1\over{|y|}}\Biggr)^{n}.$$
Now, as in the previous section, we want to show that $M_n^w(\Delta)=0$.\par
We have already shown in the previous sections that 
$w(y)=\Biggl(\ln\,\ln {1\over{|y|}}\Biggr)^{n}\in A_n(\Omega)$. 
Thus, all we need to show according to Lemma 2.3 is that the improper integral
$$\int_{\Bbb B^n(0,{1\over {e^e}})} |y|^{(1-n)\, {n\over {n-1}}}\,
\Biggl(\ln\,\ln {1\over{|y|}}\Biggr)^{{n\over{1-n}}}\,dy$$
diverges. Taking spherical coordinates in $\Bbb R^n$ the above integral becomes
$$\int_0^{{1\over {e^e}}} r^{-1}\,\Biggl(\ln\,\ln {1\over{r}}\Biggr)^{{n\over{1-n}}}\,dr.$$
Let us use the substitution $u=\ln {1\over r}$ then $du=-{1\over r}\,dr$, and the integral becomes
$$\int_{e}^\infty (\ln u)^{{n\over{1-n}}}\, du.$$
Since we have that $\ln u\le u^\epsilon$ for any $\epsilon$ positive and $u$ large enough, and thus
$$\int_{ e}^\infty (\ln u)^{{n\over{1-n}}}\, du\ge \int_{e}^\infty (u)^{\epsilon\,{n\over{1-n}}}\, du$$
$$\lim_{a\to\infty} \Biggl[ {{u^{\epsilon\,{n\over{1-n}}+1}}\over {\epsilon\,{n\over{1-n}}+1}}\Biggr]_e^a.$$
By choosing $\epsilon$ so that $\epsilon\,{n\over{1-n}}+1>0$, the above limit is equal to infinity. Thus, 
it is enough to choose $0<\epsilon<{{n-1}\over n}$ and according to Lemma 2.3 $M_n^w(\Delta)=0$. Since 
$K(x)\in L^{n-1}_{loc}(\Omega)$, we have that $M_{n-1}^{w_1}(\Delta_1)=0$, where 
$$w_1(x)=\Biggl(\ln\,\ln {1\over{|F(x)|}}\Biggr)^{n-1}.$$
Because of the relation between weighted $p$-modulus and variational weighted $p$-capacities, 
we have that $(n-1,w_1)-\text{cap}(E)=0$. It remains to show in section $\S$8 that this implies that 
the one dimensional Hausdorff measure of $E$ is zero and that will complete the proof of our result.

\beginsection{ $\S$8.  Proof of the main result}

We want to show that the one dimensional Hausdorff measure of $E$ is zero. Let us denote the one dimensional 
Hausdorff measure of $E$ by
$$\Lambda_1(E)=\lim_{\delta\to 0} \Bigl[ \inf \lbrace \sum_i r_i\colon E\subset \bigcup 
\Bbb B^n(x_i,r_i),\,0<r_i<\delta\rbrace\Bigr],$$
where the infimum is taken over all coverings of $E$ by balls of radii less than $\delta$. By Lemma 2.2 in section 
$\S$2 we have that for the weight $w_1$
$$r^{-(n-1)}\,\int_{\Bbb B^n(x_0,r)} \Biggl(\ln\,\ln {1\over{|F(x)|}}\Biggr)^{n-1}\,dx$$
$$ \le C\,((n-1)-w_1)-\text{cap} (\Bbb B^n(x_0,r),\Bbb B^n(x_0,2r)),$$
where $C$ is a constant independent of $x_0$ and $r$.\par
Let $\Bbb B^n(0,\eta)$ and 
$\Omega_\eta=F^{-1} (\Bbb B^n(0,\eta))$. Now, consider a ring $\Bbb B^n(x_0,2r)\setminus\Bbb B^n(x_0,r)$
completely contained in $\Omega_\eta$ centered at $x_0\in E$. Observe that $|F(x)|\le \eta$ 
for all $x\in \Omega_\eta$. Therefore, we have the following inequality
$$r^{-(n-1)}\,\int_{\Bbb B^n(x_0,r)} \Biggl(\ln\,\ln {1\over{\eta}}\Biggr)^{n-1}\,dx=
C\,\Biggl(\ln\,\ln {1\over{\eta}}\Biggr)^{n-1}\, r$$
$$ \le 
r^{-(n-1)}\,\int_{\Bbb B^n(x_0,r)} \Biggl(\ln\,\ln {1\over{|F(x)|}}\Biggr)^{n-1}\,dx.$$
Now, let us define $h(x_0,r)$ as follows
$$h(x_0,r)=r^{-(n-1)}\,\int_{\Bbb B^n(x_0,r)} \Biggl(\ln\,\ln {1\over{|F(x)|}}\Biggr)^{n-1}\,dx,$$
it is immediate to see that
$$C\,\Biggl(\ln\,\ln {1\over{\eta}}\Biggr)^{n-1}\le {{h(x_0,r)}\over r}.$$
Also observe that when $\eta\to 0$ we have that $x_0\to E$ and $r\to 0$. Thus, we have that 
$\lim_{r\to 0,\,x_0\to E} {{h(x_0,r)}\over r}$, which implies that for any $\epsilon >0$ there exists a $\delta>0$ 
such that $r<\epsilon\,h(x_0,r)$ whenever $r<\delta$ and $x_0$ is close enough to $E$. Let now 
$\lbrace \Bbb B^n(x_i,r_i)\colon x_i\in E, \, 0\le r_i<\delta\rbrace$ be a covering of the set $E$. If we define by 
$$\Lambda_1^\delta(E)= \inf \lbrace \sum_i r_i\colon E\subset \bigcup 
\Bbb B^n(x_i,r_i),\,0<r_i<\delta\rbrace,$$
where without loss of generality we can assume that all the $x_i$'s are in $E$, we have that $\Lambda_1(E)=
\lim_{\delta\to 0} \Lambda_1^\delta(E).$ If we choose $\delta$ as above, we have that
$$\Lambda_1^\delta(E)\le \sum_i r_i< \sum_i \epsilon\, h(x_i,r_i)$$
$$\le \epsilon \,C\, 
\lbrace \sum_i ((n-1)-w_1)-\text{cap} (\Bbb B^n(x_i,r_i),\Bbb B^n(x_i,2r_i))\rbrace .$$
We already know by the previous section that $((n-1)-w_1)-\text{cap} (E)=0$. Hence, by the definition of the weighted 
variational capacity and using rings to cover $E$ instead of balls (observe that we can always assume 
that both rings and balls are centered at points of $E$), we have that for any $\tilde\epsilon>0$ 
we can find a covering of $E$ by rings such that $E\subset \bigcup_i (\Bbb B^n(x_i,r_i)\setminus\Bbb B^n(x_i,2r_i))$ 
and 
$$\sum_i ((n-1)-w_1)-\text{cap} (\Bbb B^n(x_i,r_i),\Bbb B^n(x_i,2r_i)) \le ((n-1),w_1)-\text{cap} (E)+\tilde\epsilon.$$
combining the above inequalities we have that
$$\Lambda_1^\delta(E)\le C\,\epsilon\,\Bigl[ ((n-1),w_1)-\text{cap} (E)+\tilde\epsilon\Bigr]$$
and since both $\epsilon$ and $\tilde \epsilon$ are arbitrary, letting $\delta\to 0$ we obtain that 
$\Lambda_1(E)=0$. In particular $F^{-1}\lbrace 0 \rbrace=E$ can not contain a segment and thus, it is titally 
disconnected. replacing $F(x)$ by $F(x)-b$ in the above argument it follows that for any $b$, 
$F^{-1} \lbrace b \rbrace $ is totally disconnected. The mapping $F$ is therefore an orientation 
preserving light mapping and it follows from a theorem of Titus and Young, see [TY], that the mapping $F$ is 
open and discrete.

\vskip 0.2in

Now, we will consider the more general case in which the multiplicity function $N(F,\Omega,y)$ is not necessarily 
essentially bounded (quasi-light). In this case we have the inequality for the modulus

$$M_{n-1}^{w_1}(\Delta_1)\le C\,\Biggl[\int_{F(\Omega)} (\rho(y))^n\,
\Biggl(\ln\,\ln {1\over{|y|}}\Biggr)^{n}\,N(F,\Omega,y)\,dy\Biggr]^{{{n-1}\over n}}$$
$$\Biggl[ \int_\Omega K(x)^{n-1}\,dx\Biggr]^{{1\over n}}.$$

Our goal now will be to find which condition 
is necessary on $N(F,\Omega,y)$ which still guarantees that the first factor on the right hand side of the 
above inequality goes to zero, since that will imply that $M_{n-1}^{w_1}(\Delta_1)=0$ due to the fact that 
$K(x)\in L^{n-1}_{loc}(\Omega)$. For this, let us define the following metric as in section $\S$5, 
$\rho=|\nabla \Bigl(\ln {1\over{|y|}}\Bigr)^\delta|$ with $\delta$ positive and strictly less 
than $1-{1\over n}$. 

We want to show that, somehow we can 
use these metrics to obtain an estimate of the modulus $M_{n-1}^{w_1}(\Delta_1)$. \par
It is also clear that for each $\gamma\in \Gamma(\Delta)$ we have that $\int_\gamma \rho\,ds=\infty \ge 1$. Thus, 
it is admissible for the family of curves $\Delta$ and thus so are $\rho_\epsilon=\epsilon\, 
|\nabla \Bigl(\ln {1\over{|y|}}\Bigr)^\delta|$.\par
Without loss of generallity we can assume that $F(\Omega)\subset \Bbb B^n(0,{1\over 2})$. 
Let us compute now 
$$\int_{\Bbb B^n(0,{1\over 2})} (\rho_\epsilon(y))^n\,
\Biggl(\ln\,\ln {1\over{|y|}}\Biggr)^{n}\,N(F,\Omega,y)\,dy.$$
Passing to spherical coordinates and considering 
$$N(F,\Omega,r)=\sup_{y\in \partial \Bbb B^n(0,r)} 
N(F,\Omega,y),$$ 
the above integral is less than or equal to
$$C\, \int_0^{1\over 2} (\rho_\epsilon(y))^n\,
\Biggl(\ln\,\ln {1\over{|y|}}\Biggr)^{n}\,N(F,\Omega,r)\,dr$$
$$=C\, \int_0^{1\over 2}\epsilon^n\,\delta^n\,
\Bigl(\ln {1\over{r}}\Bigr)^{n(\delta-1)}\,
{1\over {r}}^n\,r^{n-1}\,\Biggl(\ln\,\ln {1\over{r}}\Biggr)^{n}\,N(F,\Omega,r)\,dr$$
where $C$ is a constant independent of $\delta$, $r$, and $\epsilon$. Using the change of variable 
$u=\ln {1\over{r}}$ the above integral is transformed to
$$C\,\epsilon^n\,\delta^n\,\int_{\ln 2}^\infty
\Bigl(u\Bigr)^{n(\delta-1)}\,(\ln u)^n\, N(F,\Omega,u)\,du.$$
The fact that $0<\delta<1-{1\over n}$ implies that $n(1-\delta)>1$ thus $\Bigl(u\Bigr)^{n(\delta-1)}\,(\ln u)^n$ 
is integrable, which allows for the multiplicty function $N(F,\Omega,u)$ to be unbounded and yet the above 
improper integral to be convergent. Hence, by letting $\epsilon $ go to zero we will show that 
$M_{n-1}^{w_1}(\Delta_1)=0$ and the same argument we used above will conclude that the mapping 
$F$ is discrete and open. Hence we have proved the following result
\proclaim{Theorem 8.1} Let $F\in W^{1,n}_{loc}(\Omega;\Bbb R^n)$ be a 
nonconstant mapping whose dilatation $K(x)$ 
is in $L^{n-1}_{loc}(\Omega)$. Let $N(F,\Omega,y)$ be its multiplicity function and we define 
$N(F,\Omega,r)=\sup_{y\in \partial \Bbb B^n(0,r)} 
N(F,\Omega,y)$. Then if we have that the improper integral
$$\int_{\ln 2}^\infty
\Bigl(u\Bigr)^{n(\delta-1)}\,(\ln u)^n\, N(F,\Omega,u)\,du$$
where $u=\ln {1\over r}$ converges, then the mapping $F$ is discrete and open.
\endproclaim
Examples of unbounded $N(F,\Omega,r)$ for which the above improper integral converges are $N(F,\Omega,r)=
\Bigl(\ln\,\ln {1\over r}\Bigr)^p$ for any positive $p$ and some positive $r$'s. 
Observe that in our last result we are not 
assuming that the mapping $F$ is quasi-light. So, our last result is on the direction of Iwaniec and Sverak's 
conjecture.\par
It will be interesting to study the behavior of $N(F,\Omega,r)$ as $r$ tends to zero for mappings 
$F\in W^{1,n}_{loc}(\Omega;\Bbb R^n)$ whose dilatation $K(x)$ 
is in $L^{n-1}_{loc}(\Omega)$. If somehow we would be able to show that for any of those mappings, the 
improper integral 
$$\int_{\ln 2}^\infty
\Bigl(u\Bigr)^{n(\delta-1)}\,(\ln u)^n\, N(F,\Omega,u)\,du$$
converges where $u=\ln {1\over r}$, this will prove the full Iwaniec and Sverak's conjecture.

\remark{Remark} In [Ball] it was conjectured that if a mapping $F\in W^{1,n}_{loc}(\Omega;\Bbb R^n)$ be a 
nonconstant sense preserving mapping whose dilatation $K(x)$ is in $L^{n-1}_{loc}(\Omega)$ and if $F_0\colon 
\bar\Omega\to \Bbb R^n$ be a continuous mapping in $\bar \Omega$ and one to one in $\Omega$ such that 
$F=F_0$ on $\partial\Omega$ then the mapping is discrete and open. \par

All we need to show is that for those mappings, $N(F,\Omega,y)$ is essentially bounded by one. 
This follows from Proposition 6 in [HM] and the fact shown in [Ball] that $\text{deg}(F,\Omega,y)=1$ 
for any $y\in F(\Omega)$ and $\text{deg}(F,\Omega,y)=0$ for any $y\in \Bbb R^n\setminus F(\bar\Omega)$. 
Proposition 6 in [HM] shows that $N(F,\Omega', y)= \text{deg}(F,\Omega',y)$ for a.e. 
$y\in \Bbb R^n\setminus F(\partial \Omega')$. Thus Ball's 
conjecture follows from our last theorem in this section $\S$8.
\endremark


\Refs
\widestnumber\key{HKMJ}

\ref 
\key AD
\by Adams, D. R.
\paper  Weighted nonlinear potential theory
\jour Trans. Amer. Math Soc.
\vol 297
\yr 1986 
\pages 73--94 
\endref



\ref 
\key B 
\by Beurling, A.
\paper  Ensembles exceptionnels
\jour Acta Mathematica
\vol 72
\yr 1940 
\pages 1--13 
\endref

\ref
\key  CFMS 
\by Caffarelli, L., Fabes E., Mortola S., and Salsa, S.	
\paper	Boundary behavior on nonnegative solutions of elliptic 
operators
in divergence form
\jour	Indiana Univ.\ Math.\ J.\
\vol	30
\yr      1981
\pages	621--640
\endref


\ref 
\key  C
\by Carleson, L.
\book  Selected Problems in Exceptional Sets
\publ Van Nostrand Mathematical Studies 
\vol 13
\yr  1967
\endref

\ref 
\key  Ch
\by Chua, S.-K.
\book  Extension theorems in weighted Sobolev spaces
\publ Thesis, Rutgers University 
\yr  1990
\endref



\ref 
\key  HK
\by Heinonen, J. and Koskela, P.
\paper Quasiconformal maps in metric spaces with controlled geometry
\jour Acta Mathematica
\toappear
\endref



\ref
\key HKM
\by Heinonen, J., Kilpel{\" a}inen, T. and Martio, O.
\book Nonlinear Potential Theory 
\publ Oxford University Press 
\year 1993
\endref

\ref 
\key  K
\by Kilpel\"ainen, T.
\paper Weighted Sobolev spaces and capacity 
\jour Annales Academiae Scientiarum Fennicae
\vol 19
\yr 1994
\pages  95-113
\endref


\ref 
\key  KMV
\by Koskela, P., Manfredi, J., Villamor, E.
\paper  Regularity theory and traces of $\Cal A$-harmonic functions
\jour Tran. Amer. math. Soc.
\vol 348 (2)
\yr 1996
\pages  755-766
\endref

\ref
\key MV
\by Manfredi, J.\ and Villamor, E.
\paper  Traces of monotone Sobolev functions
\jour Jour. of Geom. Anal.
\vol 6 (3)
\yr 1996
\pages 433-444
\endref



\ref 
\key V\"a 
\by V\"ais\"al\"a, J.
\book Lectures in $n$-dimensional Quasiconformal Mappings 
\publ Springer Verlag
\vol Lecture Notes 229
\yr 1971  
\endref



\ref 
\key  Z
\by Ziemer, W.
\book Weakly Differentiable Functions  
\publ Springer Verlag 
\vol Graduate Text in Mathematics 120
\yr 1989 
\endref
\endRefs
\enddocument